\documentclass[11pt,a4paper,reqno]{amsart} 

\usepackage{tikz}
\usepackage{amsmath,bm,bbm, amssymb}
\usepackage{fullpage}
\usepackage{color}
\usepackage{todonotes}
\usepackage{ulem}
\usetikzlibrary{calc}
\usetikzlibrary{patterns}
\usetikzlibrary{arrows}
\usetikzlibrary{decorations.pathreplacing}
\usepackage[utf8]{inputenc}
\usepackage{pgfplots}
\usepackage{enumerate}

\definecolor{wiasblue}   {cmyk}{1.0, 0.60, 0, 0}
\definecolor{mlugreen}{RGB}{0,81,51}

\newcommand{\conv}[2][n]{\underset{#1\rightarrow #2}{\longrightarrow}}

\def\Z{\mathbb{Z}}
\def\E{\mathbb{E}}
\def\P{\mathbb{P}}

\def\R{\mathbb{R}}
\def\mc{\mathcal}
\def\ms{\mathsf}

\def\e{\varepsilon}

\def\t{\tau}
\def\g{\gamma}
\def\k{\kappa}
\def\ff{\infty}

\def\de{\delta}

\def\es{\varnothing}
\def\one{\mathbbmss{1}}
\def\ms{\mathsf}
\def\mc{\mathcal}

\def\a{\alpha}

\def\Efr{E_{\ms b}}

\def\Era{E_{\ms d}}
\def\Env{E^{\ms V}_n}
\def\Eno{E^{\ms O}_n}
\def\Envo{E^{\ms V, *}_n}
\def\Enoo{E^{\ms O, *}_n}
\def\Cnv{2}
\def\Cno{2}

\def\Cb{\mathbf C} 

\definecolor{halfgray}{gray}{0.55}
\definecolor{webgreen}{rgb}{0,.5,0}
\definecolor{webbrown}{rgb}{.6,0,0}
\definecolor{Maroon}{cmyk}{0, 0.87, 0.68, 0.32}
\definecolor{RoyalBlue}{cmyk}{1, 0.50, 0, 0}
\definecolor{Black}{cmyk}{0, 0, 0, 0}
\definecolor{pinkish}{RGB}{255, 192, 203}
\definecolor{orange}{rgb}{216, 64, 0}

\def\been{\begin{enumerate}}
\def\enen{\end{enumerate}}
\def\im{\item}
	\def\maxx{\ms{max}}
	\def\minn{\ms{min}}
	\def\nmax{N_{n, \ms{\check C}}^{\maxx}}
	\def\nmaxo{N_{1, \ms{\check C}}^{\maxx}}
	\def\pmax{\Phi_{n,  \ms{\check C}}^{\maxx}}
	\def\pmaxk{\Phi_{n,  \ms{\check C}}^{\maxx}}
	\def\pmaxo{\Phi_{1, \ms{\check C}}^{\maxx}}
	\def\nmaxl{N_{s_{n, +}}^{\maxx}}
	\def\wl{W_{s_{n,+}}}

	\def\fmaxd{\bar F_{\ms{\check C}}}
	\def\fmax{\fmaxd}
	\def\fmaxv{\bar F_{d - 1, \ms{VR}}}
	\def\umaxm{\ell_n^{\maxx, -}}
	\def\umaxp{\ell_n^{\maxx, +}}
	\def\fmaxp{\t_n^-}
	\def\fmaxm{\t_n^+}
	\def\nminv{N_{n, k, \ms{VR}}^{\minn}}
	\def\nminc{N_{n, k, \ms{\check C}}^{\minn}}
	\def\pminv{\Phi_{n, k, \ms{VR}}^{\minn}}
	
	\def\pminck{\Phi_{n, k, \ms{\check C}}^{\minn}}
	\def\pmincc{\Phi_{n, k, \ms{\check C}}^{\minn}}
	\def\nmin{\nminv}
	\def\nminl{N_{s_n^+, k}^{\minn,\ms{VR}}}
	\def\nmincl{N_{s_n^+, 1}^{\minn, \ms{\check C}}}

	\def\fmin{F_{k, \ms{VR}}}
	\def\fminc{F_{k, \ms{\check C}}}
	\def\fmincv{F_{k, \ms{\check C}/\ms{VR}}}
	\def\fminv{\fmin}
	\def\uminm{\ell_n^{\minn, -}}
	\def\fminm{\t_n^{\minn, -}}
	\def\uminp{\ell_n^{\minn, +}}
	\def\fminp{\t_n^{\minn, +}}
	\def\bet{\begin{theorem}}
	\def\ent{\end{theorem}}
\def\bec{\begin{corollary}}
	\def\enc{\end{corollary}}
	\def\bep{\begin{proof}}
	\def\enp{\end{proof}}
	\def\f{\frac}
	
	\def\PP{\mc P}
	\def\tb{t_{\ms b}}
	\def\td{t_{\ms d}}
	\def\d{\mathrm d}
	\def\beal{\begin{align}}
	\def\beals{\begin{align*}}
	\def\enal{\end{align}}
		\def\sm{\setminus}
		\def\bel{\begin{lemma}}
		\def\enl{\end{lemma}}
		\def\rcv{r_c^{\ms V}}
		\def\rco{r_c^{\ms O}}
		
		\def\r{\rho}

		\renewcommand\le{\leqslant}
		\renewcommand\ge{\geqslant}

\newcommand{\hide}[1]{}

		\def\co{\colon}

\newtheorem{theorem}{Theorem}[section]
\newtheorem{corollary}[theorem]{Corollary}

\newtheorem{lemma}[theorem]{Lemma}
\newtheorem{proposition}[theorem]{Proposition}
\theoremstyle{definition}

\newtheorem{remark}[theorem]{Remark}

\date{\today}

\begin{document}

\title{Extremal life times of persistent loops and holes}
\author{Nicolas Chenavier}
\author{Christian Hirsch}
\address[Nicolas Chenavier]{    Universit\'e du Littoral C{\^o}te d'Opale,          EA 2597, LMPA, 50 rue Ferdinand Buisson, 62228 Calais, France}
	                 \email{nicolas.chenavier@univ-littoral.fr}
\address[Christian Hirsch]{University of Groningen, Bernoulli Institute, Nijenborgh 9, 9747 AG Groningen, The Netherlands.}
\email{c.p.hirsch@rug.nl}

\begin{abstract}
		Persistent homology captures the appearances and disappearances of topological features such as loops and holes when growing disks centered at a Poisson point process. We study extreme values for the life times of features dying in bounded components and with birth resp.~death time bounded away from the threshold for continuum percolation. First, we describe the scaling of the minimal life times for general feature dimensions, and of the maximal life times for holes in the \v Cech complex. Then, we proceed to a more refined analysis and establish Poisson approximation for large life times of holes and for small life times of loops. Finally, we also study the scaling of minimal life times in the  Vietoris-Rips setting and point to a surprising difference to the \v Cech complex. 
\end{abstract}

\maketitle
\section{Introduction}
\label{int_sec}

One of the core challenges in statistics is to provide scientifically rigorous methods to detect structure among noise. While for low-dimensional datasets, we do indeed have a versatile toolkit at our disposal, the situation is radically different in high dimensions. In this context, \emph{topological data analysis (TDA)} is emerging as an exciting novel approach. The simple and effective idea is to leverage invariants from algebraic topology to extract insights from data. While initially TDA started off as an eccentric mathematical idea, it is now applied in a wide variety of disciplines, ranging from astronomy and biology to finance and materials science \cite{gidea2018topological,rien,elastic,hiraoka}.

A core component of TDA is \emph{persistent homology}. Here, we start growing balls centered at the points of a dataset. At the beginning, all points are in isolated components, but as the balls grow, new topological features such as loops or holes may appear. As the radii increase, these features disappear again, until finally everything is contained in a single connected component without loops or holes. Through this mechanism, we can associate a birth time and a death time with each such feature.

\begin{figure}[!htpb]
	\centering
	\input{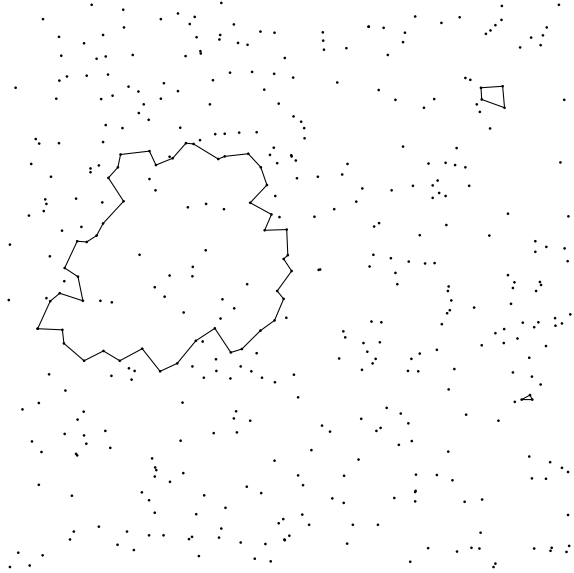}
	\caption{\v Cech-features of a maximal (left), typical (right top) and minimal (right bottom) life times.}
	\label{cyc_fig}
\end{figure}

Despite the spectacular successes of TDA in the application domains, major parts of it are still lacking a rigorous statistical foundation.  For instance, when analyzing persistent homology, practitioners often look for features living for exceptionally long periods of time, and then draw conclusions if they do occur. However, how can we decide whether the observed long life times come from genuinely interesting phenomena and are not a mere incarnation of chance? 

To answer this question, we need to understand how long life times behave under the null model of complete spatial randomness, i.e., a Poisson point process. Figure \ref{cyc_fig} illustrates cycles of a maximal, typical and minimal life times in a bounded sampling window. In this regard, a breakthrough was achieved in \cite{bobr}, where the order of the longest life time occurring in a large window was determined on a multiplicative scale. Despite the significance of this discovery, it is still difficult to build a rigorous statistical test, knowing only the order of the life time. A further interesting cross-connection is that for 0-dimensional features, the longest life time corresponds to the longest edge in the minimum spanning tree whose asymptotic was established in \cite{longest}.

In the present paper, we move one step closer to statistical applications and establish Poisson approximation results for extremal life times of loops and holes in large sampling windows. More precisely, in Theorem \ref{max_cech_thm} we determine the scaling of the maximal life time of $(d - 1)$-dimensional features, i.e., holes, in the \v Cech filtration. Moreover, we show through the Chen-Stein method that after suitable normalization, the locations of exceedances over suitable thresholds converge to a homogeneous Poisson process. Thus, the investigation is in a similar spirit as the extreme value analysis of geometric characteristics of inradii in tessellations \cite{extr1,chenavier}. However, after proceeding from the general set-up of the Chen-Stein method to more concrete tasks, the geometric analyses of the two problems are entirely different.

Another vigorous research stream in spatial extremes concerns limit results for exceptionally small structures \cite{extr3}. Our second and third main results, concern this setting. More precisely, we determine the minimal life time of features in any dimension both in the Vietoris-Rips (VR) and in the \v Cech (\v C) filtration. Our analysis reveals a striking difference in the scaling of minimal features between the VR- and the \v C-filtration. Often the VR-filtration is chosen as an approximation to the \v C-filtration, when the latter becomes prohibitively time-consuming to compute. In this respect, our result provides a clear example where this approximation fails on a fundamental level.

	The rest of the manuscript is organized as follows. First, in Section \ref{mod_sec}, we fix the notation and state the main results. Section \ref{out_sec} illustrates in broad strokes how the Chen-Stein method can be used to tackle the assertion. This involves two critical steps, namely identifying the correct scaling and excluding multiple exceedances, which will be established in Sections \ref{un_sec} and \ref{mult_sec}, respectively.	To ease the flow of reading, we have postponed overly technical volume computations into an appendix.

\section{Model and main results}
\label{mod_sec}

In this section, we describe extremal life times in TDA. First, we review briefly the fundamental concepts of the VR- and \v C-filtrations, referring the reader to \cite{chazal} for an excellent textbook treatment on this topic.

Both the VR- as well as the \v C-filtration describe collections of families of simplices in $\R^d$. A $k$-simplex $\{x_0, \dots, x_k\} \subset \R^d$ is contained in the VR-filtration at a level $r > 0$ if the distance between any pair of points is smaller than $2r$. It is included in the \v C-filtration if $\cap_{i \le k} B_r(x_i) \ne \es$, where $B_r(x_i)$ denotes the Euclidean ball with radius $r$ centered at $x_i$.

Having explained these basic building blocks, we now introduce more precisely the concept behind birth- and death times when dealing with Poisson point processes in the entire Euclidean space. For finite complexes  there is a precise algorithmic description of when a birth occurs and when a death occurs. Each time we add a simplex to the filtration it either increases the homology in the dimension of the feature or it decreases the homology in one dimension lower \cite[Algorithm 9]{chazal}. In the first case, we add a \emph{positive simplex}, which corresponds to the \emph{birth} of a feature.  In the second case, we add a \emph{negative simplex} corresponding to the \emph{death} of a feature. Defining the life time $L_{i, k}$ of the $i$th feature as the difference between the death and the birth time, this algorithm gives rise to a collection of birth-and-death times $\{L_{i, k}\}_{i \ge 1}$. Assigning to the $i$th feature the center of gravity $Z_{i, k}$ of the simplex at its death time, we may therefore consider features as a marked point process.

%
%
Let $\PP := \{X_i\}_{i \ge 1}$ be a Poisson point process in $\R^d$ with intensity 1 and $W_n := [-n^{1/d}/2, n^{1/d}/2]$ be a sampling window of volume $n$. In the \v C-filtration, persistent homology tracks the topological changes of the Boolean model
$$O_r(\PP) := \bigcup_{i \ge 1} B_r(X_i)$$
for increasing values of the radius $r$.
 For instance, 1-dimensional features are loops, and  $(d - 1)$-dimensional features correspond to {holes}, i.e., to bounded connected components in the vacant set 
$ \R^d \setminus O_r(\PP).$

%
%
Before moving to the general setting, we first describe life times for the specific case of $(d - 1)$-dimensional features. When the growing radii cause an existing hole to split into two at a radius $r$, we say that $r$ is the \emph{birth time} of the new hole. Moreover, the death time of a hole is the smallest radius $r  >  0$ when it is covered completely by the Boolean model $O_r(\PP)$. So far, this definition is still ambiguous: when a new hole is born, it is not clear, which of the newly created vacant components should determine its death time. Therefore, it is a common convention in TDA to impose that features born first die last. 

%
%
Now, we enumerate the holes as $\{H_i^*\}_{i \ge 1}$ and associate with each such hole its life time $L_i^*> 0$  and the point $Z_i^* \in \R^d$ as the last point that is covered at the death time. In this way, $\{(Z_i^*, L_i^* )\}_{i \ge 1}$ defines a stationary marked point process of holes together with their life times.

%
%
Long-range interactions coming from percolation effects make it difficult to study this marked point process in full generality. For holes born in large vacant components, the seemingly innocent convention that features born first die last may induce long-range interactions in order to decide which birth times correspond to which death times. This makes analyzing extremes prohibitively challenging: we are not aware of any form of Poisson approximation reflecting critical phenomena in percolation contexts. 

%
%
To state a clean Poisson-approximation theorem that applies to the bulk of the considered holes, we restrict $\{H_i^*\}_{i \ge 1}$ to a smaller family of holes $\{H_i\}_{i \ge 1}$ as follows. Let $\rcv  >  0$ denote the critical threshold of vacant continuum percolation and fix a time window $ [\rcv - \e_c, \rcv + \e_c]$ around $\rcv$. Then, we consider only holes with birth time outside $[\rcv - \e_c, \rcv + \e_c]$. As made precise in Section \ref{out_sec} below, the striking advantage of this restriction is that the bounded connected components are of logarithmic size \cite{cPerc}.

%
%
In the first main result, we describe the limiting behavior of exceedances of life times in the \v C-filtration centered in the sampling window $W_n$. {To do it, let $\fmax$ be the survival function of the typical life time of a hole in the \v C-filtration, i.e.,
\begin{align}
\label{thresh_max_eq}
	\fmax(\ell) = \E\big[\#\{  Z_{i}\in W_1: {L_{i} \ge \ell}\}\big].
\end{align}
We show in Appendix \ref{inv_sec} that the function $\fmax: (0, \ff) \to (0, 1)$ is invertible.
In what follows, we deal with the asymptotic behavior of the rescaled \emph{point process of exceedances}
\[\pmaxk := \big\{(n^{-1/d}Z_{i}, n\fmax(L_{i}))\big\}_{Z_{i} \in W_n} \subset W_1 \times (0, \ff).\] 
We consider this point process because it can capture both the location as well as the magnitude of large life times of features.

%
%
Our first main result shows that the maximal life time of holes in the \v Cech-filtration is of order 
$$\g_n := (\k_d^{-1} \log n)^{1/d},$$
independent of the feature dimension $1 \le k \le d - 1$, where $\k_d$ is the volume of the unit ball in $\R^d$.  Moreover, $\pmax$ is approximately a homogeneous Poisson point process $\Phi$ with intensity 1 {in $W_1 \times (0, \ff)$}. We write $\Rightarrow$ for convergence in distribution. We now motivate the quantity $\fmax^{-1}(\t/n)$ appearing the next quantity. Considering quantities of the form $\bar F^{-1}(\t/n)$, as we do in our theorems, is classical in Extreme Value Theory. Indeed, if $\{Y_i\}_i$ denotes a family of iid random variables with survival function $\bar F$, then $\P(\max \{Y_1, \dots, Y_n\} \le u_n) \to \exp(-\t)$ if and only if $n \bar F(u_n) \to \t$. 

\bet[Maximal life times, \v Cech filtration]
\label{max_cech_thm}
For every $\t > 0$,
$$\lim_{n \to \ff} \frac{\fmax^{-1}(\t/n)}{\g_n} = 1 \quad \text{and}\quad
\pmax \Rightarrow \Phi.$$
\ent
As a corollary, we obtain the scaling of the maximal life time.
\bec[Scaling of maximal life times -- \v C-filtration]
\label{max_cech_cor}
It holds that 
$$\g_n^{-1}\max_{Z_i \in W_n }L_i\xrightarrow{n \to \ff} 1$$
in probability.
\enc
\bep
First, $\max_{Z_i \in W_n}L_i \ge (1 + \e)\g_n$ means that $\{n\fmax(L_i)\}_{Z_i \in W_n}\cap [0,  n\fmax((1 + \e)\g_n)] \ne \es$. Now, we claim that $n\fmax((1 + \e)\g_n)$ tends to 0. Indeed, otherwise there exists $\t > 0$ such that $\fmax^{-1}(\t/n) \ge (1 + \e)\g_n$ for infinitely many $n$, contradicting the scaling in Theorem \ref{max_cech_thm}. Thus, we conclude by the point-process convergence in Theorem \ref{max_cech_thm}.

Similarly, $\max_{Z_i \in W_n}L_i \le (1 - \e)\g_n$ means that $\{n\fmax(L_i)\}_{Z_i \in W_n}\cap [0,  n\fmax((1 - \e)\g_n)] = \es$. We claim that $n\fmax((1 - \e)\g_n)$ tends to $\ff$. Otherwise, there exists $\t < \ff$ such that $\fmax^{-1}(\t/n) \le (1 - \e)\g_n$ for infinitely many $n$, again contradicting Theorem \ref{max_cech_thm}, and we conclude by the convergence in Theorem \ref{max_cech_thm}.
\enp

\begin{remark}
	We stress that the point-process convergence in 		Theorem \ref{max_cech_thm} 
		gives much more than just the scaling of the maximal life time: it allows us to describe the joint distribution of life times exceeding different thresholds. A downside in our definition of $\pmax$ is that the normalization of life times of the form $n\fmax(L_i)$ is rather implicit. This is because even approximating the survival function of the life time of a typical hole sufficiently accurately is challenging. If this could be achieved, then one could proceed as in \cite{extr1} to obtain a more explicit transformation. In that sense, the leading-order scaling derived in Theorem \ref{max_cech_thm} is a promising first step, but a finer analysis would be highly desirable.
	\end{remark}
	\begin{remark}
		Theorem \ref{max_cech_thm} reveals a striking difference to the multiplicative persistence approach of measuring life times in \cite{bobr}, where the scaling  is  $(\log n / \log \log n)^{1/(d - 1)} \gg \g_n$.
		 That work considers the ratio between the death time and the life time, whereas we work with the difference.
The two scalings are different as the multiplicative approach is sensitive with respect to cycles with a very small birth time. Depending on the application context this may or may not be desirable. Hence, it is valuable to know the scaling of both variants and especially how they differ.
	\end{remark}
%
%

%
%
Next, we proceed to minimal life times, where we also derive Poisson approximation results. Similarly to the setting of holes, also for general feature dimensions, the long-range correlations close to the critical radius $\rco$ of continuum percolation of the Boolean model would pose a formidable obstacle for proving Poisson approximation. Therefore, we henceforth fix again a time window $ [\rco - \e_c, \rco + \e_c]$ and consider only features with death time outside $[\rco - \e_c, \rco + \e_c]$ and that die in a bounded connected component. Furthermore, for instance by taking the center point of the $(k + 1)$-simplex whose insertion kills the $k$-dimensional feature, we attach to any such feature with life time $L_{i, k}$  a center point $Z_{i, k} \in \R^d$. Hence, we again obtain a stationary marked point process $\{(Z_{i, k}, L_{i, k})\}_{i \ge 1}$.

Fixing a feature dimension $1 \le k \le d - 1$, we define the rescaled \emph{point process of undershoots}
\[\pminck := \big\{(n^{-1/d}Z_{i, k}, n\fminc(L_{i, k}))\big\}_{Z_{i, k} \in W_n} \subset W_1 \times (0, \ff),\] 
where
\begin{align}
	\fminc(\ell) = \E\big[\#\{  Z_{i, k}\in W_1: {L_{i, k} \le \ell}\}\big].
\end{align}
denotes the distribution function of the typical life time in the \v C-filtration. In the next result, we show that the minimal life time of $k$-dimensional features is of order $n^{-2}$, independent of $k$. Moreover, the rescaled point process of undershoots is approximately a homogeneous Poisson point process with intensity 1 in $W_1 \times (0, \ff)$. However, again due to involved geometrical configurations occurring in the \v C-filtration, we prove the point-process convergence only for a specific feature dimension, namely $k = 1$ (loops).
%
%
\bet[Minimal life times, \v Cech filtration]
	\label{min_cech_thm}
	Let $1 \le k \le d - 1$. Then, for every $\t > 0$,
	\been
	\im[(i)]
	 $$\lim_{n \to \ff} \frac{\log \fminc^{-1}(\t/n)}{\log n} = -2\quad \text{and}\quad (\log n)^{-1}\min_{Z_{i, k} \in W_n }\log L_{i, k}\Rightarrow -2.$$
\im[(ii)]	 Moreover, 
$\pmincc \Rightarrow \Phi$ for $k = 1$.
\enen
\ent

%
%
The construction of $\pminck$ outlined above does not only work for the \v C- but also for the VR-filtration $\pminv$. The difference in filtration manifests itself prominently in a different scaling for the threshold. 
\bet[Minimal life times -- VR-filtration]
	\label{min_vr_thm}
	Let $1 \le k \le d - 1$. Then, for every $\t > 0$,
	$$\lim_{n \to \ff} \frac{\log \fminv^{-1}(\t/n)}{\log n} = -1 \quad \text{and}\quad \pminv \Rightarrow \Phi.$$
\ent
As a corollary, we obtain again the scaling of the minimal life time.
\begin{remark}
	A priori, it is not clear whether typical life times in the VR-filtration should be shorter or longer than those in the \v C-filtration. However, the specific construction in the proof of Theorem \ref{max_cech_thm} suggests it seems plausible that asymptotically the corresponding life-time scale $\fmaxv^{-1}(\t/n)$ should smaller than $\fmax^{-1}(\t/n)$ by a constant factor. Loosely speaking, the characteristic of persistent features is a very late death time rather than a very early birth time, and features are dying later in the \v C- than in the VR-filtration.
	\end{remark}

\section{The Chen-Stein method}
\label{out_sec}
The proofs of the Poisson approximation in Theorems \ref{max_cech_thm}--\ref{min_vr_thm} rely fundamentally on the Chen-Stein method in the form of \cite[Theorem 1]{arratia}. To make the presentation self-contained, we sketch now how to tune the general machinery to the setting of extremal life times. The delicate problem-specific conditions are then verified in Sections \ref{un_sec} and \ref{mult_sec}. 

More precisely, proceeding as in \cite{chenavier}, we discretize the problem to the Chen-Stein framework in the form of \cite[Theorem 1]{arratia}. For the convenience of the reader, we outline the most important steps for the maximal life times, noting that the arguments carry over seamlessly to the setting of minimal life times. 

%
%
The first central ingredient is to extract a dependency graph with bounded degrees. It is at this point that the percolation conditions described in Section \ref{mod_sec} enter the stage. Loosely speaking, we introduce a discretization that is coarse enough so that the life times of holes in non-adjacent blocks are independent. Hence, the size of the discretization is essentially determined by the diameter of the largest hole that we expect to see in the window $W_n$. However, due to the convention that features born first should die last, the actual block size that needs to be considered is a bit larger. More precisely, writing $S(x, r)$ for the diameter of the connected component of the vacant set $\R^d \sm O_r(\PP)$ containing $x \in \R^d$, we set 
$$S_i := \max_{r \not \in [\rcv - \e_c, \rcv + \e_c]:\, S(Z_{i}, r) < \ff} S(Z_{i}, r)$$
as the diameter of the largest bounded hole containing $Z_{i}$. 

Now, in Lemma \ref{perc_lem} in the appendix, we show that when defining the event 
\begin{align*}
	        \Env := \{&\text{for every radius $r \not \in [\rcv - \e_c, \rcv + \e_c]$ all vacant bounded connected components} \\
		        &\text{centered in $W_n$ have diameter at most $(\log n)^{\Cnv}$} \},
\end{align*}
then, for all large $n$,
\begin{align}
	        \label{enpEq}
		\P(\Env) \ge 1 - n^{-4d}.
\end{align}
Notice that, on the event $\Env$ for each hole with $Z_{i} \in W_n$, we have $L_i < (\log n)^2$.

%
%
Therefore, we partition the sampling window $W_n$ into $\mc N_n \in O(n/{(\log n)^{2d}})$  sub-boxes $\{\Cb_j\}_{j\leq \mc N_n}$ of side length 
$$s_n := 2(\log n)^{\Cnv}.$$  
 To prove that the point process of exceedances converges to a Poisson point process, we have to estimate the number of \textit{exceedances}. To do it, we first give some notation. 
For each $j\leq \mc N_n$ and $B\subset W_1$, we write
\[M_{j,n}(B):=    \max_{Z_{i} \in \Cb_j\cap (n^{1/d}B)}\big\{L_{i}:\,  S_i \le (\log n)^{\Cnv}\big\}\]
for the maximum life time of a hole in $\Cb_j$ with diameter-bound at most $(\log n)^{\Cnv}$.
Finally, we denote by 
\begin{align*}
	N_n(\tau,B) := \pmax\big((0,\tau)\times B\big)\quad\text{ and  }\quad
	N'_n(\tau,B) := \#\{j\le \mc N_n:\, n\fmaxd(M_{j,n}(B)) \le \t\}.
\end{align*}
the number of exceedances w.r.t.~$B$ (resp.~the number of sub-boxes intersected by $n^{1/d}B$ for which there exists at least one exceedance with diameter-bound lower than $(\log n)^{\Cnv}$). Setting $s_{n, +} = (3s_n)^d$, the main difficulty is to exclude multiple exceedances in small boxes of size $s_{n, +}$.

\begin{proposition}[No multiple exceedances]
	\label{mult_prop}
	Let $\tau >0$.
	\begin{enumerate}
		\item[(i)] Let $q_n:=\P\big(N_n(\tau, W_{s_{n, +}})\ge 2\big)$. Then $q_n\in O(n^{-1}).$
		\item[(ii)] The same property holds for the minimal life time in the VR-filtration with $1 \le k \le d - 1$, and for the minimal lifetime in the \v C-filtration with $k=1$.
		\end{enumerate}
\end{proposition}
We are now prepared to prove the Poisson approximation in Theorem \ref{max_cech_thm}.

%
%
\begin{proof}[Theorem \ref{max_cech_thm}, Poisson approximation]Let $\t>0$ and $B\subset W_1$ Borel.
{According to Kallenberg's theorem (see e.g. Theorem A.1, p.309 in \cite{LLR}), it is sufficient to prove that, 
\begin{equation}
\label{eq:kallenberg1}
	\E[N_n(\tau,B)] = \tau |B|
\end{equation}
and
\begin{equation}
\label{eq:kallenberg2}
	\P(N_n(\tau,B) = 0) \conv[n]{\infty} e^{-\tau|B|}.
\end{equation}
 
Equation \eqref{eq:kallenberg1} holds since 
\begin{align*}
	\E[N_n(\tau,B)] & = \E\big[ \#\{Z_{i}\in n^{1/d}B:\,  n\fmaxd(L_{i})<\tau \} \big]\\
& = n|B| \E\big[ \#\{  Z_{i}\in W_1: L_{i}>(\fmaxd)^{-1}\left( {\tau}/ n  \right)   \}  \big]\\
& = n|B|\fmaxd\left((\fmaxd)^{-1}\left( {\tau} /n  \right) \right)\\
& = \tau|B|.
\end{align*}

To deal with \eqref{eq:kallenberg2}, we first notice that
\begin{align*}
	\P\big(N_n(\tau,B) \neq N'_n(\tau,B)\big)  &\le  \P((\Env)^c) + \sum_{j \le \mc N_n} \P\big( \#\{Z_{i}\in \Cb_j:\, n\fmaxd(L_{i}) > \t\}\ge 2\big) \\
	&\le \P((\Env)^c)+\mc N_n q_n ,
\end{align*}
where the last uses that $\PP$ is stationary. Therefore, thanks to Equation \eqref{enpEq} and Proposition \ref{mult_prop}, we have $\P\big(N_n(\tau,B) \neq N'_n(\tau,B)\big)\rightarrow 0$. Thus, it is enough to prove that $\P(N'_n(\tau,B)=0)$ converges to $e^{-\tau|B|}$. We will actually prove that $N'_n(\tau,B)$ converges to a Poisson random variable with parameter $\tau|B|$. To do it, we denote by $F_n$ a Poisson random variable with parameter $\E[N'_n(\tau,B)]$. With the same computations as \cite{chenavier}, we can easily prove (thanks to Proposition \ref{mult_prop}) that 
\[\E[F_n] = \E[N'_n(\tau,B)] \underset{n\rightarrow\infty}{\sim}\E[N_n(\tau,B)]=\tau|B|.\]
Thus, it is sufficient to approximate $N'_n(B)$ by $F_n$. The main idea to do it is to apply the Chen-Stein method. Indeed, thanks to a result due to Arratia, Goldstein and Gordon (see  \cite[Theorem 1]{arratia}), and because each box $\Cb_j$ is adjacent to at most $3^d - 1$ other boxes, we have 
\begin{align}
	\label{arr_eq}
	\sup_{\mc A \subset \Z_{\ge 0}}|\P(N'_n(\tau,B) \in \mc A) - \P(F_n \in \mc A)|\le 3^d \left(2\t p_n^{(1)} + \frac{ 2n}{(2(\log n)^{\Cnv})^d}\cdot p_n^{(2)}\right),
\end{align}
where 
	$$p_n^{(1)}:= \P(n\fmaxd(M_{j, n}(W_1)) \le \t) \quad \text{ and } \quad p_n^{(2)}:=\P(n\fmaxd(M_{j_1, n}(W_1) \wedge  M_{j_2, n}(W_1)) \le \t)$$
denote the exceedance	probabilities in a single box $\Cb_j$ or in two adjacent boxes $\Cb_{j_1}$ and $\Cb_{j_2}$. By stationarity and symmetry, notice that $p_n^{(1)}$ and $p_n^{(2)}$ do not depend on $j, j_1$ or $j_2$. The first term of the right-hand side converges to 0 since, by stationarity,
	\[p_n^{(1)} \le \E\big[\#\{ Z_{i} \in \Cb_j:\, L_{i} > \fmaxd^{-1}(\t/n)\}\big] = \frac{|\Cb_j|\t}n.\]
 The second one converges to 0 according to Proposition \ref{mult_prop}. This concludes the proof of Theorem \ref{max_cech_thm}.
 }
\end{proof}
In a similar way, we obtain the Poisson approximation asserted in Theorem \ref{min_cech_thm}  and Theorem \ref{min_vr_thm}, relying on a correspondingly defined $\Eno$ instead of $\Env$.
}

\section{Scalings of $\fmax$, $\fminv$ and $\fminc$}
\label{un_sec}
In the present section, we derive the scalings of $\fmax$, $\fminv$ and $\fminc$ asserted in Theorems \ref{max_cech_thm}--\ref{min_vr_thm}.

%
%
To begin with, we deal with $\fmax$. As observed in \cite[Remark 3.2]{bobr}, the upper bound follows from the results in \cite{kahle}: according to the scaling of the coverage radius in \cite{hall}, above $\g_n$ the union of balls covers $W_n$, which is contractible. In particular, by the nerve theorem \cite[Theorem 10.7]{bjorn}, the \v C-complex has the same Betti numbers as the union of balls, so that any cycle is trivial. The lower bound resorts to a construction inspired by \cite[Lemma 5.1]{bobr}. \\[1ex]
\noindent\emph{Scaling of $\fmax$ (proof of Theorem \ref{max_cech_thm}, scaling).}
\phantom{a}\\[1ex]
Henceforth, we write 
$$\nmax(\t):= \pmaxk(W_1 \times [0, \t])$$
for the number of features in $W_n$ with life time at least $\fmax^{-1}(\t/n)$.

%
%
\noindent	{\boldmath $\fmax^{-1}(\t/n)/\g_n\le 1 + \e.$\unboldmath}
	For $\e >0$ 
	write $\umaxp := (1 + \e)\g_n$ and $\fmaxp := n\fmax(\umaxp)$. To prove that $\lim_{n \to \ff} \fmaxp= 0$, we have to prove that
	$$\E\big[\nmaxo(\fmaxp/n)\big] \in o(n^{-1}).$$
	In order to bound the expected number of holes with death time at least $\umaxp$ that are centered in $W_1$, let 
$R := \min\{|x|\co{x \in \PP} \}$
		be the distance of the closest point of $\PP$ to the origin. In particular, the death time of any hole centered in $W_1$ is at most $R + d$, as $W_1$ is covered completely at the latest at this time. Next, there can be at most one feature with life time exceeding $\umaxp$  and centered in $W_1$. Indeed, if there were at least two such features $H_1, H_2$, then at least one of them is born after time $R - d$, as otherwise, $H_1$ and $H_2$ would still be the same hole at the corresponding birth time. Combining this with the previous observation on the death time leads to the contradiction that such a hole would have life time at most $2d$. 	Hence, 
			$$\E\big[\nmaxo( \fmaxp/n)\big]= \P\big(\pmaxo \cap (W_1 \times [0, \fmaxp/n]) \ne \es\big) \le \P\big(\PP \cap B_{\umaxp - d}(o)  = \es\big).$$
				Invoking the void probability of the Poisson process $\PP$, the right-hand side is at most 
					$$\exp\big(-\k_d (\umaxp - d)^d\big) =  \exp\big(-(1 + \e)^d \log(n) (1 - d/\umaxp)^d\big) \in o(n^{-1}).$$
					 \\[1ex]

%
%
%
	%

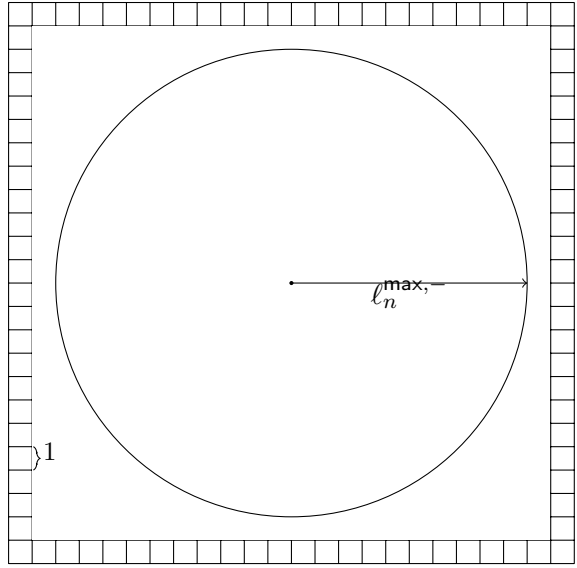
\begin{figure}[!htpb]
		\centering
	\begin{tikzpicture}[scale=1.55]

	\draw[step=.2] (-2.4, -2.4) grid (2.4,2.4);
	\fill[white] (-2.195, -2.195) rectangle (2.195, 2.195);


	\fill (0,0) circle (.5pt);
	\draw (0,0) circle (2cm);
	\draw[->] (0,0)--(2,0);

	\draw[decorate, decoration={brace}] (-2.19,-1.4)--(-2.19,-1.6);

	\coordinate[label={$\umaxm$}] (A) at (1, -0.30);
	\coordinate[label=\small{1}] (A) at (-2.05, -1.6);
\end{tikzpicture}
		\caption{Template for maximal life time feature.}
		\label{max_fig}
	\end{figure}

	%
	%
	\noindent	{\boldmath $\fmax^{-1}(\t/n)/\g_n\ge 1 - \e.$\unboldmath}
	 With $\umaxm:= (1 - \e)\g_n$ and $\fmaxm := n\fmax(\umaxm)$, we claim that $\lim_{n \to \ff}\fmaxm= \ff$, i.e.,
	$$\lim_{n \to \ff}n\E\big[\nmaxo( \fmaxm/n)\big] = \lim_{n \to \ff}\E\big[\nmax(  \fmaxm)\big] = \ff.$$
	As indicated at the beginning of this section, we modify a construction of \cite[Lemma 5.1]{bobr}. Loosely speaking, we cover a discretized annulus of radius $\umaxm$ by Poisson points, see Figure \ref{max_fig}. 

	%
	%

	More precisely, let 
	$$\Efr := \big\{\PP \cap W(v) \ne \es \text{ for all $v \in \Z^d$ with $|v|_\ff = 2\lceil  \g_n \rceil$}\big\}$$
denote the event that each unit block in a $d_\infty$-annulus of radius $2\lceil \g_n \rceil$ contains at least one Poisson point. We also write $m_n = m_{n, \e}$ for the number of these blocks and note that $m_n \in O((\log n)^{(d - 1)/d})$. 
	Moreover, we let 
$$\Era := \big\{\PP \cap B_{\umaxm}(o) = \es\big\}$$
	denote the  event that there are no Poisson points in the ball with radius ${\umaxm}$. 
	Then,  under the event $\Efr \cap \Era$ we find a hole with life time at least $\umaxm$. 

	Setting $\e' := 1 - (1 - \e)^d$ the independence property of the Poisson process gives that
	$$\P(\Efr \cap \Era) = \P(\Efr) \P(\Era) \ge e^{-c m_n}n^{-1 + \e'} $$
	for some $c  > 0$.
	Finally, we can take the template described by $\Efr \cap \Era$ and shift it to a different location in the window $W_n$. Since those configurations are of logarithmic extent, we can arrange at least $a_n := n^{1 - \e'/2}$ of them disjointly in $W_n$. 
	In particular, since $m_n \in O((\log n)^{(d - 1)/d})$, 
	$$\lim_{n \to \ff}\E\big[\nmax( \fmaxm)\big] \ge \lim_{n \to \ff} a_ne^{-c  m_n}n^{-1 + \e'}  = \ff.$$
	
%
%
Next, we analyze the scaling of $\fminc$. Since we now deal with short rather than long life times, the proof structure is converse to that for $\fmax$. That is, in the lower bound we analyze configurations leading to small life times, whereas the upper bound relies on a specific construction.
	We also write
$$\nminc(\t):= \pminck(W_1 \times [0, \t])$$
for the number of features in $W_n$ with life time at most $\fminc^{-1}(\t/n)$.

Henceforth, let 
\begin{align}
	\label{filt_eq}
	d_{0\cdots k}:= d(x_0, \dots, x_k):= \inf\{r > 0\co \bigcap_{x \in \{x_0, \dots, x_k\}}B_r(x) \ne \es\}
\end{align}
denote the \emph{filtration time} of the simplex $\{x_0, \dots, x_k\}$, i.e., the first time it appears in the \v C-filtration. In particular,  $d_{ij} = |x_i - x_j|/2$.

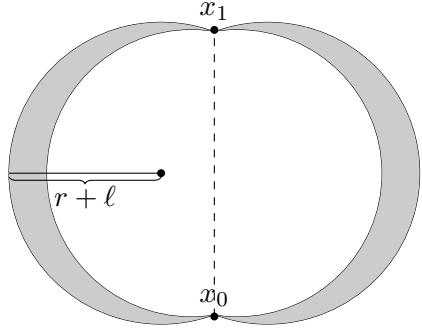
\begin{figure}[!htpb]
	\centering
	\begin{tikzpicture}[scale=1]

	\draw (-2.2, 0) circle (2cm);
	\fill[black!20!white] (-2.2, 0) circle (2cm);
	\draw (-1.8, 0) circle (1.9cm);

	\draw (-0.8, 0) circle (2cm);
	\fill[black!20!white] (-0.8, 0) circle (2cm);
	\draw (-1.2, 0) circle (1.9cm);
	\fill[white] (-1.2, 0) circle (1.9cm);
	\fill[white] (-1.8, 0) circle (1.9cm);

	\draw[dashed] (-1.5,-1.9)--(-1.5,1.9);
	\draw (-2.2, 0)--(-4.2, 0);
	\draw[decorate, decoration={brace}] (-2.2, -.05)--(-4.2, -.05);

\fill (-1.5, -1.9) circle (1.5pt);
\fill (-1.5, 1.9) circle (1.5pt);
\fill (-2.2, 0) circle (1.5pt);



\coordinate[label={$x_0$}] (A) at (-1.5, -1.9);
\coordinate[label={$x_1$}] (A) at (-1.5, 1.9);
	\coordinate[label={$r + \ell$}] (A) at (-3.2, -.6);
\end{tikzpicture}
	\caption{The crescent set $D_{r, \ell}(x_0, x_1)$.}
	\label{cresc_fig}
\end{figure}

Moreover, if $x_0, \dots, x_{k - 1} \in \R^d$ are affinely independent and $r > d(x_0, \dots, x_{k - 1})$, then in every $k$-dimensional hyperplane containing $x_0, \dots, x_{k - 1}$, there are precisely two $k$-dimensional balls with radius $r$ containing $x_0, \dots, x_{k - 1}$ on their boundary, and we let $D_r(x_0, \dots, x_{k - 1})$ denote the union of all such balls. Then, for $r, \ell > 0$, we introduce the crescent
$$D_{r, \ell}(x_0, \dots, x_{k - 1}) := D_{r + \ell}(x_0, \dots, x_{k - 1}) \setminus D_r(x_0, \dots, x_{k - 1}),$$
illustrated in Figure \ref{cresc_fig} for $k = 2$.
 having a $k$-simplex $x_0, \dots, x_k$ with filtration time $s$ means that their exist $i_0, \dots, i_m \le k$ and an $m$-dimensional ball $B_s(P)$ centered at some $P \in \R^d$ such that $x_0, \dots, x_k \subset B_s(P)$ and $x_{i_0}, \dots, x_{i_m} \in \partial B_s$
 . In particular, $s \in [r, r+ \ell]$ means that  $x_{i_m} \in D_{r, \ell}(x_{i_0}, \dots, x_{i_{m - 1}})$.

\noindent\emph{Scaling of $\fminc$ (proof of Theorem \ref{min_cech_thm}, scaling).}
\phantom{a}\\
%
%
	{\boldmath $\fminc^{-1}(\t/n) \ge n^{-2 - \e}\;\text{\bf and }\;\P\big(\min_{Z_{i, k} \in W_n }L_{i, k} \le n^{-2 -  \e} \big) \to 0.$\unboldmath}
	It suffices to prove the first assertion.
We write $\uminm := n^{-2 -\e}$ and $\fminm := n\fminc(\uminm)$ and claim that $\lim_{n \to \ff} \fminm= 0$, i.e.,
	$$\lim_{n \to \ff}\E\big[\nminc(\fminm)\big] = 0.$$
When a \v C-feature with life time $\ell > 0$ is born at time $r = d_{0\cdots k}$, then there exist Poisson points $X_0, \dots, X_k$ such that $\cup_{i \le k}B_r(X_i)$ covers the simplex spanned by these points. If a feature dies at time $\td = d_{0\cdots k} + \ell$, then there is $m \le k $ and Poisson points $X_0', \dots, X_{m + 1}'$ such that the associated simplex is covered at time $d_{0\cdots k} + \ell$. After reordering and choosing $m$ to be the smallest such value, there exists $j \le  m$ such that $X_i' = X_i$ for $i \le j$ and $X_i' \ne X_i$ for $i \ge j + 1$. Then, the previous observations can be succinctly summarized as
$$X'_{m + 1} \in D_{d_{0\cdots k},  \ell}(X_0, \dots, X_j, X_{j + 1}', \dots,  X_m').$$

Now, we invoke the Slivnyak-Mecke formula from \cite[Theorem 4.4]{poisBook}. It shows that the expected number of Poisson points $X_0, \dots, X_j, X_{j + 1}', \dots,  X_{m + 1}'$ satisfying 
$$X'_{m + 1} \in D_{d_{0\cdots k},  \ell}(X_0, \dots, X_j, X_{j + 1}', \dots,  X_m')$$
is at most 
\begin{align*}
	&\int_{W_{2n} \times B_{s_n}(x_0)^{j } \times B_{s_n}(x_0)^{m - j} }\int_{D_{d_{0\cdots k}, \ell}(x_0, \dots, x_m)  }1\d x_{m + 1}' \d (x_0,  \dots, x_{j  },x_{j + 1}',  \dots, x_{m  }')\\
	&\quad= 	\int_{W_{2n} \times B_{s_n}(x_0)^{j } \times B_{s_n}(x_0)^{m - j} }\hspace{-.5cm}{}\big|D_{d_{0\cdots k}, \ell}(x_0, \dots, x_m)  \big| \d (x_0,  \dots, x_{j  },x_{j + 1}',  \dots, x_{m  }') .
\end{align*}
Now, Lemma \ref{cresc_lem} in the appendix shows that the volume of the integrand is bounded above by 
$cs_n^d\sqrt \ell,$ so that 
\begin{align*}
	\int_{W_{2n} \times B_{s_n}(x_0)^{j } \times B_{s_n}(x_0)^{m - j} }\hspace{-.5cm}{}\big|D_{d_{0\cdots k}, \ell}(x_0, \dots, x_m)  \big| \d (x_0,  \dots, x_{j  },x_{j + 1}',  \dots, x_{m  }') 
	\le c \k_d^{2m - j}s_n^{(2m - j + 1)d}\int_{W_{2n}}\hspace{-.3cm} \sqrt \ell \d x_0,
\end{align*}
and for $\ell \le \uminm$ the latter expression tends to 0 as $n \to \ff$. \\[1ex]

%
%
\noindent	{\boldmath $\fminc^{-1}(\t/n) \le n^{-2 + \e}\;\text{\bf and }\;\P\big(\min_{Z_{i, k} \in W_n }L_{i, k} \ge n^{-2 +  \e} \big) \to 0$.\unboldmath}
Next, writing $\uminp := n^{-2 +\e}$ and $\fminp := n\fminc(\uminp)$, we claim that $\lim_{n \to \ff} \fminp= \ff$, i.e.,
	$$\lim_{n \to \ff}\E\big[\nmin(\fminp)\big] = \ff.$$
The proof idea is to construct features with a short life time by relying on perturbations from a deterministic template. We first define suitable templates of constant size, which in the final argument will then be copy-pasted throughout the entire window $W_n$.
	To lower the technical barriers, we give the proof for $k = d - 1$. The arguments for the lower-dimensional settings are very similar.  

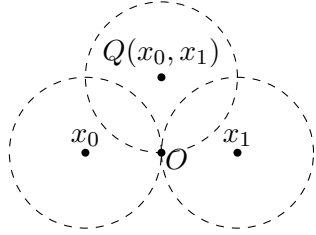
\begin{figure}[!htpb]
		\centering
	\begin{tikzpicture}[scale=1]

\fill (-1, 0) circle (1.5pt);
\fill (1, 0) circle (1.5pt);
\fill (0, 1) circle (1.5pt);
\fill (0, 0) circle (1.5pt);

	\draw[dashed] (-1, 0) circle (1cm);
	\draw[dashed] (1, 0) circle (1cm);
	\draw[dashed] (0, 1) circle (1cm);

\coordinate[label={$x_0$}] (A) at (-1, -0.05);
\coordinate[label={$x_1$}] (A) at (1, -0.05);
\coordinate[label={$O$}] (A) at (0.2, -0.35);
	\coordinate[label={$Q(x_0, x_1)$}] (A) at (0, 1);
\end{tikzpicture}
		\caption{Template for \v C-feature of minimal life time in dimension $d = 2$.}
		\label{fig_cech}
	\end{figure}

First, for affinely independent points $x_0, \dots, x_{d - 1} \in \R^d$ let $Q(x_0, \dots, x_{d - 1}) \in \R^d$ denote one of the two possible points projecting onto the center $O$ of the $(d - 1)$-dimensional circumball of $x_0, \dots, x_{d - 1}$, and whose distance from that center is given by the radius of that sphere. Figure \ref{fig_cech} illustrates this configuration in dimension $d = 2$. 

	To fix a distinguished template, we consider the unit vectors $P_i := e_{i + 1} \in \R^d$, $0 \le i \le d - 1$, forming a regular $(d - 1)$-simplex of side length $\sqrt 2$ with circumradius $\r_{d - 1} = \sqrt{(d - 1) / d}$. 
We claim that at time $\r_{d - 1}$ the faces of the simplex $\{P_0, \dots, P_{d - 1}, Q\}$ are all covered.
	Since $\r_{d - 1}$ is the coverage radius of a regular simplex, it suffices to show by symmetry that $|P_0 - Q| \le \sqrt 2$. 
	Indeed, remark that $OP_0Q$ is an isosceles rectangular triangle and write 
	$|P_0-Q|^2=2\r_{d-1}^2={2(d-1)}/d$, which gives $|P_0-Q|=\sqrt{2}\sqrt{(d-1)/d} < \sqrt 2$. 

	Next, we consider random perturbations of this template of a small magnitude $\de \in (0, 1/2)$. More precisely, 
	we define the event 
$$E := \big\{\PP(B_\de(P_i)) =  1 \text{ for all $i \le d - 1$}\big\} \cap\big\{\PP(B_{4d}(o)) = d  + 1\big\}$$ 
that $\PP$ has precisely one point in the $\de$-neighborhood of each of the points in the template in dimension $k$ and no other points in the ball $B_{4d}(o)$. Conditioned on $E$, we write $P_i'$ for the Poisson point contained in $B_\de(P_i)$, which is then uniformly distributed at random in $B_\de(P_i)$. Moreover, under the event $E$ the Poisson process $\PP$ contains precisely one further point inside $B_{4d}(o)$, which we denote by $Q'$. Hence, writing $d_{0'\cdots(d - 1)'} = \inf\{r > 0\co \bigcap_{x \in \{P_0', \dots, P_{d - 1}'\}}B_r(x) \ne \es\}$,  conditioned on $E$, the probability that this feature has life time at most $\ell$  is at least 
\begin{align*}
	&\P\big(Q' \in D_{d_{0'\cdots (d - 1)'}, \ell}(P_0', \dots, P_{d - 1}') \cap B_\de\big(Q(P'_0, \dots, P'_{d - 1})\big)\big) 	\\	
	&= \f1{\k_d((4d)^d - d\de^d)} \E\Big[\big|D_{d_{0'\cdots (d - 1)'}, \ell}(P_0', \dots, P_{d - 1}') \cap B_\de\big(Q(P'_0, \dots, P'_{d - 1})\big)\big|\Big].
\end{align*}
Then, Lemma \ref{cresc_lem} in the appendix shows that the volume on the right-hand side is at least $c\sqrt \ell$ for a suitable $c > 0$.  
Now, for a suitable $c' > 0 $, we can fit at least $c'n$ such templates into the window $W_n$. Hence, for $\ell = \uminp$,
	$$\lim_{n \to \ff}\E\big[\nminc(\fminp)\big] \ge \lim_{n \to \ff} cc'n \sqrt{\uminp}  = \ff,$$
	as asserted. To prove that $\P\big(\min_{Z_{i, k} \in W_n }L_{i, k} \ge n^{-2 +  \e} \big) \to 0$, we note that by the independence of the Poisson point process, the probability that none of the shifted templates yields a feature with life time at most $\uminp$ is at most $(1 - c\sqrt{\uminp})^{c'n}$, which tends to 0, since $n\sqrt{\uminp} \to \ff$.

%
%
Finally, we move to the scaling of $\fminv$. 
Henceforth, 
	$$B_{r, \ell}(x) := B_{r + \ell}(x) \setminus B_r(x)$$
	denotes the annulus with center $x \in \R^d$, inner radius $r$ and thickness $\ell$.\\[1ex]

\noindent\emph{Scaling of $\fminv$ (proof of Theorem \ref{min_vr_thm}, scaling).}
\phantom{a}\\
%
%
	{\boldmath $\fminv^{-1}(\t/n) \ge n^{-1 - \e}.$\unboldmath}
	Writing $\uminm := n^{-1 -\e}$ and  $\fminm := n\fminv(\uminm)$, we now assert that $\lim_{n \to \ff} \fminm= 0$, i.e.,
	$$\lim_{n \to \ff}\E\big[\nmin(\fminm)\big] = 0.$$
	First, we may assume that $\Eno$ occurs. Indeed, by Cauchy-Schwarz,
	$$\E\big[\nmin(\fminm)\one\{(\Eno)^c\}\big] \le \E\big[ \PP (W_{2n})^d\one\{(\Eno)^c\}\big] \le \big(\E\big[ \PP( W_{2n})^{2d}\big]\P((\Eno)^c)\big)^{1/2}.$$ 
	Since $\PP$ is a Poisson point process, the expectation on the right-hand side is of order $O(n^{2d})$, so that it remains to invoke a property which is similar to \eqref{enpEq}.

	%
	%
	For the remaining argument, recall that a simplex belongs to the VR-filtration at time $t$ if the distance between any pair of its vertices is at most $2t$. In particular, when a feature with life time $\ell$ is born at time $\tb$, there exist $X_1, X_2 \in \PP$ with $d_{12} = \tb$, and when it dies at time $\td = \tb + \ell$, there are Poisson points $X_3, X_4$ with $d_{34} = \tb + \ell$. In other words, $X_4 \in B_{2d_{12}, 2\ell}(X_3)$. Moreover, since the event $\Eno$ occurs whp, we may assume that $X_1 \in W_{2n}$ and that $\max_{2 \le i \le 4}d_{1i} \le s_n$. 

	%
	%
	Here, $\{X_1, \dots, X_4\}$ consists of at least three different points. First, if they are all pairwise distinct, then the Slivnyak-Mecke formula  bounds the expected number of undershoots as 
	\begin{align*}
		\int_{W_{2n}}\int_{B_{s_n}(x_1)^2}\big|B_{2d_{12}, 2\uminm}(x_3)  \big| \d( x_2, x_3) \d x_1
		&= \k_d^2s_n^d\int_{W_{2n}}\int_{B_{s_n}(x_1)}\big((2d_{12} + 2\uminm)^d - (2d_{12})^d\big) \d x_2 \d x_1\\
		&\le 4^dd\k_d^3s_n^{3d}\int_{W_{2n}}\uminm  \d x_1\\
		&= 8^dd\k_d^3s_n^{3d} n^{-\e},
	\end{align*}
	which tends to 0 as $n \to \ff$. 
	Second, if say $X_2 = X_3$, then we proceed similarly to obtain the bound
\begin{align*}
	\int_{W_{2n}}\int_{B_{s_n}(x_1)}\big|B_{2d_{12}, 2\uminm}(x_2)  \big| \d x_2 \d x_1
	&= \k_d\int_{W_{2n}}\int_{B_{s_n}(x_1)}\big((2d_{12} + 2\uminm)^d - (2d_{12})^d\big) \d x_2 \d x_1\\
	&\le 4^dd\k_d^2s_n^{2d}\int_{W_{2n}}\uminm  \d x_1\\
		&= 8^dd\k_d^2s_n^{2d} n^{-\e},
\end{align*}
which again tends to 0 as $n \to \ff$. \\[1ex]

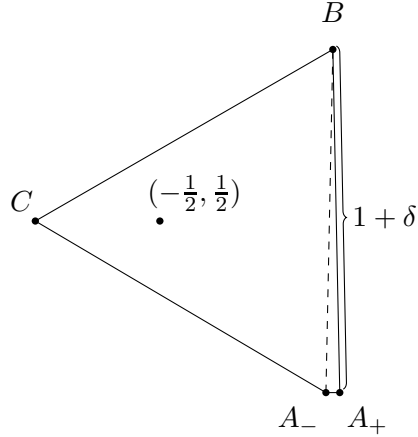
\begin{figure}[!htpb]
	\centering
\begin{tikzpicture}[scale=4.55]

\fill (-.02, 0) circle (.3pt);
\fill (.02, 0) circle (.3pt);
\fill (0, 1) circle (.3pt);
\fill (-.5, .5) circle (.3pt);
\fill (-.86, .5) circle (.3pt);

\draw (-.02, 0)--(.02, 0)--(0, 1)--(-.86,.5)--cycle;
	\draw[decorate, decoration = {brace}] (0.01, 1.01)--(.03, 0.01);
\draw[dashed] (-.02, 0)--(0,1);

\coordinate[label={$A_-$}] (A) at (-.1, -0.15);
\coordinate[label={$A_+$}] (A) at (.1, -0.15);
\coordinate[label={$B$}] (A) at (0, 1.05);
	\coordinate[label={$(-\tfrac12, \tfrac12)$}] (A) at (-.4, .5);
\coordinate[label={$C$}] (A) at (-.9, .5);
\coordinate[label={$1 + \de$}] (A) at (0.15, .45);
\end{tikzpicture}
		\caption{Template for VR-feature of minimal life time.}
		\label{min_vr_fig}
\end{figure}

%
%
\noindent	{\boldmath $\fminv^{-1}(\t/n) \le n^{-1 + \e}.$\unboldmath}
Next, writing $\uminp := n^{-1 +\e}$ and $\fminp := n\fminv(\uminp)$, we claim that  $\lim_{n \to \ff} \fminp= \ff$, i.e.,
	$$\lim_{n \to \ff}\E\big[\nmin(\fminp)\big] = \ff.$$
To that end, we construct a feature with a short life time for each feature dimension $1 \le k \le d - 1$. 
The features depend on a deterministic template illustrated in Figure \ref{min_vr_fig}.
%
%

We fix a small value $\de = \de(k, d) > 0$. Then, we set 
$$A_\pm := \pm\sqrt{(1 + \de)^2 - 1} e_1, \quad B := e_2,$$ 
so that $|A_- - B| = |A_+ - B| = 1 + \de$. We also fix a point $C := (\eta, 1/2)$, where $\eta < 0$ is chosen, so that $|C - A_+| > 1 + \de > |C - B|$.

%
%
To convey the idea, we sketch how this template gives rise to the desired feature for $k = 1$ before moving to higher $k$. If we consider the VR-complex on $\{A_-, A_+, B, C\}$ at level $1 + \de$, then $A_-A_+BC$ is a loop. Moreover, removing the edge $A_-B$ and the attached higher simplices, the complex does not contain any triangles as by construction $|C - A_+| > 1 + \de$. In particular, the loop has positive life time. However, after adding the edge $A_-B$, the loop becomes the boundary of the triangles $A_-A_+B$ and $BCA_-$. 

%
%
To generalize the construction to higher $k$, we introduce additional points 
$\{P_{1,\pm}, \dots, P_{k - 1,\pm} \}$
to the complex as follows. First, set 
$$P_{i, \pm} := \frac 12(e_2 - e_1) \pm (0.5 + 4\de) e_{i + 2}.$$
In particular, the $P_{i, \pm}$ are at distance at most $1 + \de$ from $A_-, A_+, B$ and $C$.  
Now, we generalize the above consideration for the VR-complex at level $|A_+ - B|$ but with the edge $\{A_-, B\}$ removed. In particular, the complex contains all $k$-simplices of the form 
$$\{P_{1,\e_1}, P_{2,\e_2}, \dots, P_{{k - 1},\e_{k - 1}}, P, P'\},$$
where $\e_i \in \{-, +\}$ and $\{P, P'\}$ is one of $\{A_-, A_+\}$, $\{A_+, B\}$, $\{B, C\}$ or $\{C, A_-\}$. Then, these simplices form a cycle. Indeed, when removing the vertex $P_{i,\e_i}$, then the resulting face is also contained in the $k$-simplex with $P_{i,\e_i}$ replaced by $P_{i,-\e_i}$. On the other hand, if for instance $\{P, P'\} = \{A_-, A_+\}$ and if we remove the point $A_-$, then we find the corresponding face also in the simplex with $\{P, P'\} = \{A_+, B\}$. However, this complex does not contain any $(k + 1)$-simplices. Indeed, $|P_{i,-} - P_{i,+}| > 1 + \de$ are at distance larger than $1 + \de$, and for any triple from $\{A_-, A_+, B, C\}$ at least one edge is also not in the complex. This situation changes drastically if $A_-B$ belongs to the complex. Then, the boundaries of the $(k + 1)$-simplices
$$\{P_{1,\e_1}, P_{2,\e_2}, \dots, P_{{k - 1},\e_{k - 1}}, A_-, A_+, B\} \quad \text{ and }\quad \{P_{1,\e_1}, P_{2,\e_2}, \dots, P_{{k - 1},\e_{k - 1}}, B, C, A_-\}$$
yield the cycle constructed before.

%
%
It remains to connect this template to random features induced by $\PP$. Similarly as in the \v Cech case, set $\de_1 = \de/8$ and consider the event 
$$E' := \big\{\PP(B_{\de_1}(P)) = 1 \text{ for all $P \in \{A_-, A_+, B, C\} \cup \{P_{i,\pm}\}_{i \le k - 1}$}\big\} \cap\big\{\PP(B_3(o)) = 2k + 2\big\}$$ 
that $\PP$ has precisely one point in the $\de_1$-neighborhood of each of the template points and no other points in the ball $B_3(o)$. Conditioned on $E'$, we write $A_-', A_+', B', C' $ and $(P_{i, \pm})'$ for the Poisson points lying in the corresponding neighborhoods of the template points. In particular, under $E'$, the primed points are distributed uniformly at random in the corresponding $\de_1$-neighborhoods. Moreover, the construction of the template is sufficiently robust with respect to $\de_1$-perturbations to show that the adjacency findings for the VR-complex continue to hold upon replacing the template points by the perturbed points. 

%
%
Conditioned on $E'$, the probability that this feature has a short life time is at least 
\begin{align*}
	&\P(0 < |A_-' - B'| - |A_+' - B'| < \uminp\,|\,E') 	\\	
	&= (\k_d\de_1^d)^{-3} \Big|\Big\{(A_\pm',  B')\co \max_{P \in \{A_\pm, B\}}|P' - P| \le \de_1 
	\text{ and } 0 < |A_-' - B'| - |A_+' - B'| < \uminp\Big\}\Big|.
\end{align*}
Now, we set 
$u := A_+ - B'$ { and } $v := A_-' - B',$
and for $\a \in [0, 1]$ and $\eta \in \R^d$ 
$$A_{+,\a, \eta} := B' + (|v| - \a \uminp)(u + \eta)/|u + \eta|,$$
so that
$|A_{+, \a, \eta} - B'| = |v| - \a \uminp$.
In particular,
$|A_-' - B'| - |A_{+, \a, \eta} - B'| = \a \uminp < \uminp.$

Hence, writing $\de_2 := \de_1/8$,  it suffices to show that 
$$|A_{+, \a, \eta} - A_+| \le \de_1$$
for every $A_-' \in B_{\de_2}(A_-)$, $B' \in B_{\de_2}(o)$, $\a \in [0, 1]$ and $\eta \in B_{\de_2}(o)$, because then every such point $A_{+, \a, \eta}$ is an admissible choice for $A_+'$. 
To that end, we leverage the bound
\begin{align*}
	|A_{+,\a, \eta} - A_+| &\le \Big|- u + |v|\frac u{|u|} + |v|\Big(\frac{u + \eta}{|u + \eta|} - \frac u{|u|}\Big)\Big| + \uminp\\
	&\le \big||u| - |v|\big| + |v|\Big|\frac{u + \eta}{|u + \eta|} - \frac u{|u|}\Big| + \uminp.
\end{align*}
The first expression is at most $3\de_2$ and the second one at most $4\de_2$, so that taken together, the right-hand side is indeed at most $\de_1$.

Now, note that 
$$\big|\{A_{+, \a, \eta}\co \a \in [0, 1],\, \eta \in B_{\de_2}(o)\}\big| \ge \uminp \k_{d - 1}\de_2^{d - 1},$$
so that the probability is much larger than $n^{-1}$, and consequently that the expectation tends to infinity.

\section{Proof of Proposition \ref{mult_prop}}
\label{mult_sec}
In this section, we prove that there are no multiple exceedances in $W_{s_{n, +}}$.

%
%

\bep[Proof of Proposition \ref{mult_prop} (i)]
By the scaling derived in Section \ref{un_sec}, it suffices to prove that
\begin{align*}
	\P\big(\nmaxl(s_{n, +}\fmaxd(\umaxm)) \ge 2\big) \in O(n^{-1}).
\end{align*}
%
%
By the scaling in Section \ref{un_sec}, both features $H_i$ and $H_j$ die at time at most $\umaxp$, so that their birth times are at most $2\e \g_n$. 
In particular, for the holes $H_i$ and $H_j$ to be contained in different connected components at the birth time of the later hole, say $H_i$, we have $d(Z_i, \partial H_i)\ge (1-\e)\g_n$ and $d(Z_j, \partial H_j)\ge (1-3\e)\g_n$. Since $d(Z_i,Z_j)\ge d(Z_i, \partial H_i)+d(Z_j, \partial H_j)$ at the birth time of $H_i$, this implies that 
 $ |v - v'| \ge (2 - 5\e)\g_n\ge \g_n, $ for $\e$ small enough, where  $W(v) := v + W_1, W(v') := v' + W_1$ are the lattice boxes containing $Z_i$ and $Z_j$, respectively.
Therefore,
$$\big|B_{\umaxm}(v) \cup B_{\umaxm}(v')\big| \ge \k_d^{-1}\log(n) |B_{1 - \e}(o) \cup B_{1 - \e}(e_1)|.$$

%
%

Then,
\begin{align*}
	&\P\big(\nmaxl(s_{n, +}\fmaxd(\umaxm)) \ge 2\big) \\
	&\quad\le\sum_{\substack{v, v' \in W_{s_{n, +}} \cap \Z^d\\ |v - v'| \ge \g_n}}  \P\big(L_{i} \wedge L_{j} \ge \umaxm \text{ for some } Z_i \in W(v),\,Z_j \in W(v')  \big)\\
	&\quad\le \sum_{\substack{v, v' \in W_{s_{n, +}} \cap \Z^d\\ |v - v'| \ge \g_n}}  \P\Big( \big( B_{\umaxm}(v) \cup B_{\umaxm}(v')  \big)\cap \PP=\es  \Big)\\
& \quad \le |\wl^2|\exp\big(-\k_d^{-1}\log(n) |B_{1 - \e}(o) \cup B_{1 - \e}(e_1)|\big).
\end{align*}

The last term is in $O(n^{-1})$ provided that $\e > 0$ is sufficiently small.
\enp

%
%

Next, we move to the minimum in the \v C-filtration.
\bep[Proof of Proposition \ref{mult_prop} (ii); \v C]
Again,  it suffices to show that
\beal
	\label{mult_min2_eq}
	\P(\nmincl(\uminp) \ge 2) \in O(n^{-1}),
	\end{align}
where we now  set $\uminp := n^{-2 + \e}$.
	%
	%
	We need to understand the geometric implications of finding features $Z_i, Z_j \in \wl$ with life times shorter than $\uminp$. We recall from Section \ref{un_sec} that one such undershoot gives Poisson points $\{X_0, \dots, X_4\}$ such that the triangle $X_2X_3X_4$ is covered by the union of balls centered at $X_i$, $2\le i\le 4$. First, we observe that if at least one of $X_0, X_1$ is not contained in $\{X_2, X_3, X_4\}$, then \eqref{mult_min2_eq} holds even without taking the second feature $Z_j$ into account. Indeed, if for instance $X_1 = X_2$, but $X_0 \not \in \{X_2, X_3, X_4\}$, then the Slivnyak-Mecke formula yields the bound
	$$\int_{\wl^3}|B_{2d_{234} - 2\uminp, 2\uminp}(x_2)| \d (x_2, x_3, x_4) \in O(s_n^{4d} \uminp).$$

	%
	%
	We may hence assume to have points $\{X_0, X_1, X_2\}$ with $d_{01} < d_{012} \le d_{01} + \uminp$ for the first feature and similarly points $\{X_0', X_1', X_2'\} \ne \{X_0, X_1, X_2\}$ with $d_{0'1'} < d_{0'1'2'} \le d_{0'1'} + \uminp$. Again, several configurations are possible, where the most challenging ones correspond to the cases where $\{X_0, X_1, X_2\}$ and $\{X_0', X_1', X_2'\}$ differ in one variable. 
	
%
	%
	First, assume that $X_0' = X_0$, $X_1' = X_1$, but $X_2'\ne X_2$. Then, invoking Lemma \ref{cresc_lem} in Appendix \ref{apdx_sec} gives the volume bound
	\begin{align*}
		\Big|\big\{(x_0, x_1, x_2, x_2') \in \wl^4:\,d_{012}  ,d_{012'}\in (d_{01}, d_{01} + \uminp) \big\}\Big|
		&\le \int_{\wl^2}\big|D_{d_{01}, \uminp}(x_0, x_1)\big|^2 \d(x_0, x_1)\\
		&\le c_{\ms{cresc}}^2s_{n, +}^4\uminp.
	\end{align*}
	Noting that the last line is of order $O(n^{-1})$ concludes the proof of the first case.
	Second, assuming that $X_1' = X_1$, $X_2' = X_2$, but $X_0' \ne X_0$, Lemmas \ref{cresc_lem} and \ref{int_min2_lem} in Appendix \ref{apdx_sec} give that 
	\begin{align*}
		&\Big|\big\{(x_0, x_1, x_2, x_0') \in \wl^4:\,d_{01}<d_{012} \le d_{01} + \uminp \text{ and } d_{0'1} < d_{0'12} \le d_{0'1} + \uminp\big\}\Big|\\
		&\quad\le \int_{\wl^2}\int_{D_{d_{01}, \uminp}(x_0, x_1)}\big|\{x_0' \in \wl:\,d_{0'1}<d_{0'12} \le d_{0'1} + \uminp \}\big| \d x_2\d(x_0, x_1)\\
		&\quad\le cn^{-1/4}\int_{\wl^2}\big|D_{d_{01}, \uminp}(x_0, x_1)\big|\d(x_0, x_1)\\
		&\quad \le cc's_{n, +}^3n^{-1/4}\sqrt{\uminp},
	\end{align*}
	which is again in $O(n^{-1})$.
	Finally consider the case where $X_1' = X_2$, $X_2' = X_1$, but $X_0' \ne X_0$. Then, similarly, Lemmas \ref{cresc_lem} and \ref{int_min2_lem} in Appendix \ref{apdx_sec} give that 
	\begin{align*}
		&\Big|\big\{(x_0, x_1, x_2, x_0') \in \wl^4:\,d_{01}<d_{012} \le d_{01} + \uminp \text{ and } d_{0'2} < d_{0'21} \le d_{0'2} + \uminp\big\}\Big|\\
		&\quad\le \int_{\wl^2}\int_{D_{d_{01}, \uminp}(x_0, x_1)}\big|\{x_0' \in \wl:\,d_{0'2}<d_{0'21} \le d_{0'2} + \uminp \}\big| \d x_2\d(x_0, x_1),
	\end{align*}
	so that we now conclude as in the previous case.
\enp

%
%
Finally, we deal with the minimum in the VR-filtration.
\bep[Proof of Proposition \ref{mult_prop} (ii); VR]
Again, it suffices to prove that
\beal
	\label{mult_min_eq}
	\P(\nminl(\uminp ) \ge 2) \in O(n^{-1}), 
	\end{align}
where we now  set $\uminp  = n^{-1 + \e}$. Now, suppose two features centered at $Z_i, Z_j \in \wl$ live shorter than $\uminp$. By the VR-filtration, this means that for the feature centered in $Z_i$, there exist Poisson points $X_0, \dots, X_3 \in \wl$ such that 
$$d_{01} < d_{23} < d_{01} + \uminm.$$
Similarly also the feature centered in $Z_j$ gives rise to Poisson points $X'_0, \dots, X'_3 \in \wl$ such that 
$$d_{0'1'} < d_{2'3'} < d_{0'1'} + \uminm.$$
Again, not all points need to be distinct, but both $\{X_0, \dots, X_3\}$ and $\{X_0', \dots, X_3'\}$ consist of at least 3 elements. Moreover, $\{X_0, X_1\} \ne \{X_0', X_1'\}$ and $\{X_2, X_3\} \ne \{X_2', X_3'\}$ since the features are distinct.

We now distinguish two cases. The first one being that $\{X_0, \dots, X_3\} \ne \{X_0', \dots, X_3'\}$. Say, for instance $X_3' \not \in \{X_0, \dots ,X_3\}$. Then, we may apply the Slivnyak-Mecke formula to see that the expected number of configurations is at most 
\begin{align*}
	&\int_{W_{s_n+}}\hspace{-.4cm}\E\Big[\#\big\{X_0, \dots, X_3, X_0',\dots, X_2' \in \PP \cap \wl\colon d_{01} < d_{23} < d_{01} + \uminm\text{ and } x_3 \in B_{2d_{01}, 2\uminm}(X_2)\big\} \Big]\d x_3'\\
	&\le 2^{d + 1} ds_{n, +}^{4d} \uminm\E\Big[\#\big\{X_0, \dots, X_3 \in \PP \cap \wl \colon d_{01} < d_{23} < d_{01} + \uminm\big\}\Big]\\
	&\le 2c(2^d ds_{n, +}^{4d})^2 (\uminm)^2 \in O(n^{-1}),
\end{align*}
for a suitable $c > 0$.
It remains to deal with the case $\{X_0', \dots, X_3'\} = \{X_0, \dots, X_3\}$, which necessarily means that they consist of 4 elements. For instance, it may occur that $X_0' = X_0$, $X_1' = X_3$, $X_2' = X_2$ and $X_3' = X_1$. Then,
$d_{01} < d_{23} < d_{01} + \uminp$ { and } $d_{03}  < d_{12} < d_{03} + \uminp.$
Applying the Slivnyak-Mecke formula, the expected number of such configurations in $\wl$ is at most 
$$
\Big|\big\{(x_0, \dots, x_3) \in \wl^4:\, |d_{23} - d_{01}| \vee|d_{03} - d_{12}| \le \uminp \big\}\Big|.
$$
Lemma \ref{int_max_lem} in Appendix \ref{apdx_sec} shows that the latter expression is in $O(n^{-1})$, as asserted.
\enp

\bibliography{lit}
\bibliographystyle{abbrv}

\appendix
\section{Invertibility of $\fmin$ and $\fminc$}
\label{inv_sec}
In this section, we show that $\fmin, \fminc:\, (0, \ff) \to (0, 1)$ are invertible for every $1 \le k \le d - 1$.

%
%
\bel[Invertibility of $\fmin$ and $\fminc$]
\label{inv_lem}
Let $1 \le k \le d - 1$. Then, $\fmin, \fminc:\, (0, \ff) \to (0, 1)$ are invertible.
\enl
\bep
We prove invertibility by showing that $\fmincv$ is a strictly decreasing function, which is continuous and satisfies $\lim_{\ell \to 0}\fmincv(\ell)  = 0$ and $\lim_{\ell \to \ff} \fmincv(\ell) = 1$. \\[1ex]
\noindent{\boldmath $\lim_{\ell \to 0}\fmincv(\ell) = 0$.} First, by right-continuity $\lim_{\ell \to 0}\fmincv(\ell) = \fmincv(0)$, where $\fmincv(0)$ describes the expected number of features with life time 0 centered in $W_1$. However, by definition, the life time of a feature is always strictly positive, so that $\fmincv(0) = 0$.\\[1ex]
\noindent{\boldmath $\lim_{\ell \to \ff}\fmincv(\ell) = 1$.} Now, $\lim_{\ell \to \ff} 1 - \fmincv(\ell)$ describes the expected number of features in $W_1$ with an infinite life time. However, this number vanishes. Indeed, in the \v C-filtration the feature must have died by the time that a representing cycle is contained in a ball, since the latter is contractible. If we take a \v C-boundary chain as a witness for the death of the feature in the \v C-filtration, then this chain is also a witness for the death in the VR-filtration. Thus, with probability 1, all life times are finite.\\[1ex]
\noindent{{\bf Continuity.}} Let $\ell > 0$ be arbitrary. To show continuity, we establish that there are no features with life time exactly $\ell$. We start by proving the claim for VR filtration. In that case, writing $d_{ij} = |X_i - X_j|/2$, there would exist points $X_0, X_1, X_2, X_3 \in \PP$ such that $d_{23} = d_{01} + \ell$ with $X_3 \not \in \{X_0, X_1, X_2\}$. Then, 
\begin{align*}
	\E\Big[\#\{X_3:\, X_3 \in \partial B_{2d_{01} + 2\ell}(X_2)  \text{ for some $X_0, X_1, X_2 \in \PP$} \}\Big] 
	&= \int_{\R^d} \P\Big(x \in \partial B_{2d_{01} + 2\ell}(X_2) \text{ for some $X_0, X_1, X_2 \in \PP$} \Big) \d x\\
	&= \E\Big[ \Big|\bigcup_{X_0, X_1, X_2 \in \PP} \partial B_{2d_{01} + 2\ell}(X_2)\Big|\Big]\\
	&= 0.
\end{align*}
For the \v C-filtration, the argumentation is similar. Indeed, we recall from Section \ref{un_sec} that a \v C-feature gives Poisson points $X_0, \dots, X_k \in \PP$ and $X_0', \dots, X_{k + 1}' \in \PP$ with $X_{k + 1}' \not \in \{X_0', \dots, X_k'\}$ such that relying on the filtration times from \eqref{filt_eq}, we have $d_{0'\cdots (k + 1)'} = d_{0\cdots k} + \ell$. Hence, similarly as in the VR filtration,
\begin{align*}
	&\E\Big[\#\{X_{k+1}':\, d_{0'\cdots (k + 1)'} = d_{0\cdots k} + \ell \text{ for some $X_0, \dots, X_k, X_0', \dots, X_k' \in \PP$} \}\Big] \\
	&= \E\Big[ \Big|\bigcup_{\substack{X_0, \dots, X_k \in \PP \\ 
	X_0', \dots, X_k' \in \PP}} \partial D_{d_{0\cdots k}}(X_0', \dots, X_k')\Big|\Big]\\
	&= 0.
\end{align*}
\\[1ex]
\noindent{{\bf Strict monotonicity.}} Strict monotonicity means that for $\ell < \ell'$ there is a positive probability to observe a feature with life time in $(\ell, \ell')$. Note that scaling all points by a factor also scales the life time by that factor. Hence, it suffices to show that for some fixed $b$ and every $\e > 0$, there is a positive probability to have a feature with life time in $(b - \e, b + \e)$. 

We start with the VR-filtration. Consider a feature described by a cross-polytope $\{\pm e_i\}_{1 \le i \le k + 1}$. Then, this feature has birth time $1/\sqrt 2$ and death time $1$. Therefore its life time is $b := 1 - 1/\sqrt 2$. If we allow the vertices of the feature to be perturbed by at most $\e/2$, the life time is in the interval $(b - \e, b + \e)$. 

Finally, we deal with the \v Cech-filtration. Consider a feature described by a regular simplex $\{e_i\}_{1 \le i \le k + 1}$. Then, this feature has birth time $\r := \sqrt {(k - 1) / k }$ and death time $\r' := \sqrt{k /(k + 1)}$. Therefore, its life time is $b := \r' - \r$,  and we conclude as in the VR-filtration.
\enp

%
%
\section{Percolation properties}
\label{perc_sec}
In this section, we show that occupied or vacant components in continuum percolation whose birth time is bounded away from the critical threshold are of poly-logarithmic size with high probability. More precisely, recalling that we center components at their center of gravity, we fix $\e_c > 0$ and  put
\begin{align*}
	\Enoo := \{&\text{for every radius $r \not \in [\rco - \e_c, \rco + \e_c]$ all occupied bounded connected components} \\
	&\text{centered in $W_n$ have diameter at most $(\log n)^{\Cno}$} \}.
\end{align*}
and 
\begin{align*}
	\Envo := \{&\text{for every radius $r \not \in [\rcv - \e_c, \rcv + \e_c]$ all vacant bounded connected components} \\
	&\text{centered in $W_n$ have diameter at most $(\log n)^{\Cnv}$} \},
\end{align*}
where $\rco, \rcv$ are the critical radii for occupied and vacant continuum percolation, respectively.

We show that the events $\Enoo$ and $\Envo$ occur with high probability. For a \emph{fixed} sub- or supercritical radius $r$, this follows from classical continuum percolation theory \cite[Chapter 4]{cPerc}. However, since in the definitions of $\Enoo$ and $\Envo$ the radius may vary, we need a small discretization argument.

%
%
\bel[$\Enoo$ and $\Envo$ occur with high probability]
\label{perc_lem}
For all sufficiently large $n \ge 1$ it holds that 
$$\P(\Enoo) \ge 1 - n^{-4d} \quad \text{ and } \quad \P(\Envo) \ge 1 - n^{-4d}.$$
\enl
\bep
We explain how to proceed in the occupied setting, noting that the arguments in the vacant case are almost identical. We say that an event occurs whp if its complement occurs with probability at most $n^{-5d}$ for large $n$.
First, below $\rco$ all connected components are almost surely bounded, and the size of these components increases with the radius. Hence, we conclude from \cite[Lemma 3.3]{cPerc} that the corresponding components are of logarithmic size whp. Second, above $\rcv$ all vacant components are bounded and their size decreases with increasing radius. Since every bounded occupied component lies within a vacant component, we conclude from \cite[Lemma 4.1]{cPerc} that again the components are of logarithmic size whp.

%
%
Hence, it remains to consider $r \in (\rco + \e_c, \rcv)$, noting that this interval is non-empty only if $d > 2$.
Here, we discretize the possible radii that can occur. More precisely, we subdivide the interval $[\rco + \e_c, \rcv]$ into $K_n \in O(n^{8d})$ parts $\{r_i\}_{i \le K_n}$ such that $r_{i + 1} - r_i \le n^{-8d}$. Then, by \cite[Theorem 10.20]{penrose}, whp for each $r_i$, all bounded connected components at level $r_i$ centered in $W_n$ are of diameter at most $(\log n)^2$. Hence, it suffices to show that for any $r \in (\rco + \e_c, \rcv)$ the components at level $r$ correspond to the components at level $r_i$ for some $i \le K_n$. In other words, we claim that whp, for each $i \le K_n$ there exists at most one pair $\{X_j, X_k\} \subset W_n$ with $|X_j - X_k|/2 \in (r_i, r_{i + 1})$.

%
%
To that end, note that for fixed $i$, the Slivnyak-Mecke formula shows that the probability that there exist distinct pairs $\{X_j, X_k\} \subset W_n$ and $\{X_j', X_k'\} \subset W_n$ with $|X_j - X_k|/2 \in (r_i, r_{i + 1})$ and $|X_j' - X_k'|/2 \in (r_i, r_{i + 1})$ is of order $O(n^{2(1 - 8d)})$. Hence, by the union bound, the probability that some  $i \le K_n$ has this property is at most $O(n^{2 - 8d})$, as asserted.
\enp

%
%
\section{Volume computations}
\label{apdx_sec}
In this section, we compute volume bounds for the specific configurations occurring in the proofs in Sections \ref{un_sec} and \ref{mult_sec}. All computations rely only on findings from elementary geometry, but are still bit tedious when written out in detail.  

First, we bound the volumes of crescents. As in the proof of Theorem \ref{min_cech_thm}, for any affinely independent points $x_0, \dots, x_{d - 1} \in \R^d$ let $Q(x_0, \dots, x_{d - 1}) \in \R^d$ denote one of the two possible points projecting onto the center of the $(d - 1)$-dimensional circumsphere of $x_0, \dots, x_{d - 1}$, and whose distance from that center is given by the radius of that sphere. As in Sections \ref{un_sec} and \ref{mult_sec}, we write $d_{i_0\cdots i_k}$ for the filtration time when the simplex $\{x_{i_0}, \dots x_{i_k}\}$ appears in the \v Cech filtration.

%
%
\bel[Volume of crescents]
\label{cresc_lem}
There exists $c_{\ms{cresc}} = c_{\ms{cresc}}(d) > 0$ with the following properties. Let $r > 0$, $ \ell < 1$, $x_0, \dots, x_{d - 1} \in \R^d$ be affinely independent. Then,
$$\big|D_{r, \ell}(x_0, \dots, x_{d - 1}) \big| \le c_{\ms{cresc}}(r + \ell)^{d - 1}\sqrt \ell.$$
Moreover, if $d_{0\cdots (d - 1)} \ge 1/2$ and $\de \in (\ell, 1)$, then
$$\big|D_{d_{0\cdots (d - 1)}, \ell}(x_0, \dots, x_{d - 1})  \cap B_\de(Q(x_0, \dots, x_{d - 1}))\big|\ge c_{\ms{cresc}}^{-1} \de^{d - 1}\sqrt \ell.$$
\enl
\bep

%
%
 To ease notation, set $a := d_{0 \cdots (d - 1)}$. 
By rotating and shifting we may assume that $x_0, \dots, x_{d - 1}$ are contained in $\{0\} \times \R^{d - 1}$ and that their circumcenter is the origin. Furthermore, the set $D_{r, \ell}(x_0, \dots, x_{d - 1})$ is rotationally symmetric around the axis $\R e_1$. 

\noindent {\bf Upper bound.} Defining $r_+(b), r_-(b) > 0$ for any $b \in \R$ through 
$$D_{r, \ell}(-a e_2, ae_2) \cap (\{b\} \times \R) = \{b\} \times \big([r_-(b), r_+(b)] \cup [-r_+(b), -r_-(b)]\big),$$
we obtain by Fubini that
\begin{align*}
	\big|D_{r, \ell}(x_0, \dots, x_{d - 1})\big| &= \int_0^{2r + 2\ell} \big|B_{r_+(b)}(o) \sm B_{r_-(b)}(o)\big| \d b
	= \k_{d - 1}\int_0^{2r + 2\ell}(r_+(b)^{d - 1} - r_-(b)^{d - 1})\d b.
\end{align*}
Now, $r_+(b)^{d - 1} - r_-(b)^{d - 1}$ is at most $d(r + \ell)^{d - 2}(r_+(b) - r_-(b))$, so that 
$$\big|D_{r, \ell}(x_0, \dots, x_{d - 1})\big| \le d\k_{d - 1}(r + \ell)^{d -2}\int_0^{2r + 2\ell}(r_+(b) - r_-(b)) \d b\\
	= \tfrac d2\k_{d - 1}(r + \ell)^{d -2}|D_{r, \ell}(-a e_2, ae_2)|.$$
Hence, we have now reduced the upper bound to the special setting where $d = 2$. Here, an elementary geometric argument that is elucidated in \cite[Lemma 9.8]{bchs} concludes the proof.

\noindent{\bf Lower bound.} For a point $x \in \R^{d - 1}$ with $|x| \le \de$, we let
$$I(x) := \{b \ge 0:\, (b, x) \in D_{a, \ell}(x_0, \dots, x_{d - 1})\}$$
denote the interval consisting of all points inside $D_{a, \ell}$ projecting onto $x$.
Then, it suffices to show that there exists a constant $c > 0$ such that for each such $x$ we have $|I(x)| \ge c \sqrt \ell$. Again, after rotation, we may reduce to the two-dimensional setting, i.e., assume that $x = (0, x, 0, \dots, 0)$ for some $x \le \de$. 

\begin{figure}[!htpb]
	\centering
	\begin{tikzpicture}[scale=1]

	\draw (-1.5, 0) circle (1.9cm);

	\draw (-.8, 0) circle (2cm);
	\draw[dashed] (-.8,0)--(1.0, -1);

	\draw[dashed] (-1.5, .7)--(1.1, .7);
	\draw[dashed] (-.8, 0) -- (-.8, .7);

	\draw[dashed] (-1.5,-1.9)--(-1.5,1.9);
	\draw[decorate, decoration={brace}] (-.8, -.05)--(-1.5, -.05);
	\draw[decorate, decoration={brace}] (1.1, .65)--(.2, .65);
	\draw[very thick] (1.08, .7)--(.25, .7);
	\draw[decorate, decoration={brace}] (-1.5, .0)--(-1.5, .7);

\fill (-1.5, 0) circle (1.5pt);
\fill (-1.5, .7) circle (1.5pt);
\fill (-.8, 0) circle (1.5pt);
\fill (-1.5, -1.9) circle (1.5pt);
\fill (-1.5, 1.9) circle (1.5pt);



\coordinate[label={$-ae_2$}] (A) at (-1.5, -2.3);
\coordinate[label={$ae_2$}] (A) at (-1.5, 1.9);
\coordinate[label={$0$}] (A) at (-1.7, -.4);
\coordinate[label={$b_0$}] (A) at (-1.2, -.5);
\coordinate[label={$x$}] (A) at (-1.7, .2);
\coordinate[label={$a + \ell$}] (A) at (-.3, -.9);
\coordinate[label={$M$}] (A) at (-.7, -.4);
	\coordinate[label={$I(x)$}] (A) at (.7, .15);
\end{tikzpicture}
	\caption{Lower bound on crescent volume}
	\label{cresc_app_fig}
\end{figure}
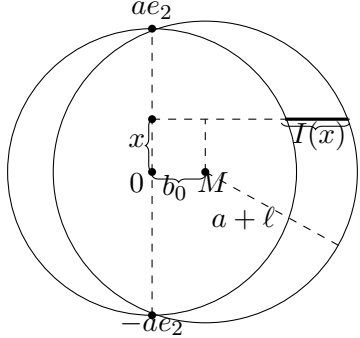

To derive the lower bound on $|I(x)|$, note that the first coordinate of one of the midpoint of the circle through $\pm a e_2$ with radius $a + \ell$ equals 
$$b_0 := \sqrt{(a + \ell)^2 - a^2} = \sqrt \ell \sqrt{2a + \ell}.$$
Thus, as illustrated in Figure \ref{cresc_app_fig},
$$|I(x)|= \sqrt{(a + \ell)^2 -x^2} + b_0 - \sqrt{a^2 -x^2} = b_0 + \ell\frac{2a + \ell}{\sqrt{(a + \ell)^2 -x^2} + \sqrt{a^2 -x^2}},$$
which is bounded below by a scalar multiple of $\sqrt \ell$ since $b_0 \ge \sqrt{2a \ell} \ge \sqrt \ell$.

\enp

Next, we bound the volume of almost parallelograms, recalling that $s_{n, +} = 3^d s_n^d$.

%
%
\bel[Almost parallelogram]
\label{int_max_lem}
If $\e > 0$ is sufficiently small, then
$$
\Big|\big\{(x_0, \dots, x_3) \in \wl^4:\, |d_{23} - d_{01}| \vee|d_{03} - d_{12}| \le n^{-1 + \e}  \big\}\Big| \in O(n^{-1}).
$$
\enl

\bep
%
%
 First, by rotating $x_3$ into the plane spanned by $x_0$, $x_1$ and $x_2$, we may reduce to the two-dimensional setting.
%
%
Second, we may assume that $\min_{i \ne j} d_{ij} \ge n^{-1/16}$. Indeed, for instance if $i  = 0$ and $j = 1$, then
$$\Big|\big\{(x_0, \dots, x_3) \in \wl^4:\, d_{01}\le n^{-1/16} \text{ and } |d_{23} - d_{01}| \le \ell_n  \big\}\Big| \in o(n^{-1}),$$
where we set $\ell_n := n^{-1 + \e}$.
 Similarly, we may also assume that the angles $\angle x_ix_jx_k$ are at least $n^{-1/16}$ for all pairwise distinct $i,j,k$. Note that we did not attempt to optimize the exponent $-1/16$.

\begin{figure}[!htpb]
	\centering
	\begin{tikzpicture}[scale=1.5]

	\begin{scope}
	\clip(-2.2,.8) rectangle (1.0,3.2);

	\draw[dashed] (-2, 1)--(.8,1);

	\begin{scope}
	\clip(.8, 1) circle (2.71cm);
	\clip(-2, 1) circle (2.41cm);
	\fill[black!30] (-2, 1) circle (2.4cm);
	\fill[black!30] (.8, 1) circle (2.7cm);
	\draw (-2, 1) circle (2.2cm);
	\end{scope}

	\fill[white] (-2, 1) circle (2.01cm);
	\fill[white] (.8, 1) circle (2.31cm);

	\begin{scope}
		\clip (-1.9, 1.8) rectangle (-1.0, 1.9);
	\draw (-2, 1) circle (1.0cm);
	\end{scope}


	\begin{scope}
	\clip(.8, 1) circle (2.71cm);
	\clip(-2, 1) circle (2.41cm);
	\clip(-2.2,1) rectangle (1.0,3.2);
\draw[] (-2, 1) circle (2.4cm);
	\draw[] (.8, 1) circle (2.7cm);

	\draw[] (-2, 1) circle (2.0cm);
	\draw[] (.8, 1) circle (2.3cm);

	\end{scope}

	        \draw[decorate, decoration = {brace}] (0.4,.97)--(.0,.97);
	        \draw[decorate, decoration = {brace}] (-1.5,.97)--(-1.9,.97);

	\draw   (0,-1) arc(63:0:-2cm);
	\draw  (.37,-.05) arc(63:85:-.9cm);
	\draw[ultra thick]   (0,-1) arc(63:81:-2cm);

	\end{scope}

	\draw[dashed] (-2.0, 1)--(-.69, 2.76);
	\draw[dashed] (-2.0, 1)--(-1.03, 2.98);

	\fill (.8, 1) circle (1.0pt);
	\fill (-2, 1) circle (1.0pt);
	\fill (-1.03, 2.98) circle (1.0pt);
	\fill (-0.69, 2.76) circle (1.0pt);

%
\coordinate[label={$x_0$}] (A) at (-2.0, .75);
\coordinate[label={$x_2$}] (A) at (.8, .75);
	\coordinate[label={$Q$}] (A) at (-1.03, 2.98);
	\coordinate[label={$P$}] (A) at (-0.59, 2.5);
\coordinate[label={$\ell_n$}] (A) at (-1.7, .65);
\coordinate[label={$\ell_n$}] (A) at (.25, .65);
\coordinate[label={$r$}] (A) at (-1.55, 1.95);
\coordinate[label={$\a$}] (A) at (-1.55, 1.65);
\end{tikzpicture}
	\hspace{2cm}
	\begin{tikzpicture}[scale=1.05]
	\clip(-4.3,-2.4) rectangle (2,2.6);

	\draw[ultra thick] (-2.2, 0) circle (2cm);
	\fill[white] (-6, -1.05) rectangle(4,4);

	\draw (-1.5, 1.9) circle (3.8cm);
	\fill[white] (-1.5, -2.0) rectangle (4.4,5);
	\fill[white] (-5.5, -1.05) rectangle (4.4,5);
	\draw (-1.5, 1.9) circle (1.3cm);
	\fill[white] (-1.35, 0) rectangle (4.4,5);
	\fill[white] (-1.95, 0) rectangle (-4.4,5);
	\draw (-2.2, 0) circle (2cm);

	\draw[decorate, decoration = {brace}] (-1.4,1.9)--(-1.4,-1.9);

	%
	%
	\draw[dashed] (-1.5, 1.9)--(-2.9,-1.9);
	\draw[decorate, decoration = {brace}](-2.95,-1.9)--(-1.55,1.9);
\fill (-1.5, 1.9) circle (1.5pt);
	\fill (-2.2, 0) circle (1.5pt);

\coordinate[label={}] (A) at (-1.5, -2.3);
\coordinate[label={$x_1$}] (A) at (-1.5, 1.9);
	\coordinate[label={{$2a - 2\ell_n$}}] (A) at (-0.75,-.2);
\coordinate[label={$\tfrac\a2$}] (A) at (-1.6,.6);
	\coordinate[label={$2a$}] (A) at (-2.55, -.1);
\end{tikzpicture}
	\caption{Intersection area of annuli (left); Arc describing those $x_0$ with $d_{01}\ge d_{012} - \ell_n$ (right)}
	\label{paral_fig}
\end{figure}
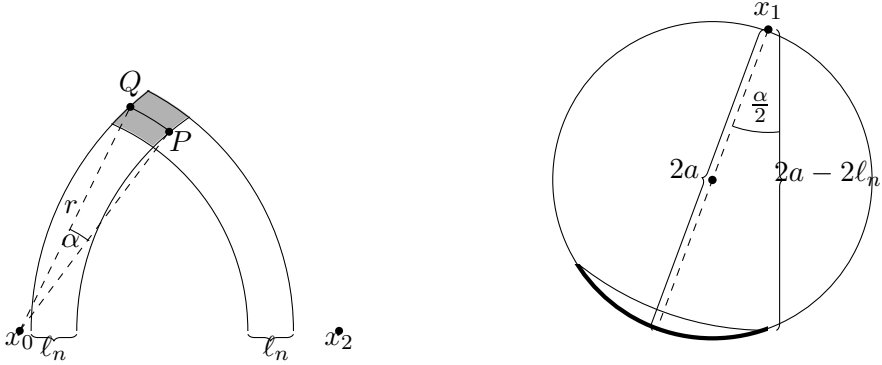

%
%
After these simplifications, it remains to bound the annuli-intersection area 
$
|B_{d_{12}, \ell_n}(x_0) \cap B_{d_{01}, \ell_n}(x_2)|,
$
as illustrated in Figure \ref{paral_fig}. To bound this quantity, we rely on the co-area formula from \cite[Chapter 3]{coarea}. More precisely, we have 
$$|B_{d_{12}, \ell_n}(x_0) \cap B_{d_{01}, \ell_n}(x_2)| = \int_{d_{01}}^{d_{01} + \ell_n} |\partial B_r(x_0) \cap B_{d_{01}, \ell_n}(x_2)|\d r.$$
%
%
Hence, if we write $P, Q$ for the endpoints of one of the two arcs $\partial B_r(x_0) \cap B_{d_{01}, \ell_n}(x_2)$ and $\a$ for the enclosed angle, it suffices to show that $\a$ is of order $O(n^{-1/2})$.

To that end, note that 
$$\cos(\angle x_2x_0P) =  \frac{d_{02}^2 + r^2 - d_{01} ^2}{2d_{02}r}\quad \text{ and }\quad\cos(\angle x_2x_0Q) = \frac{d_{02}^2 + r^2 - (d_{01} + \ell_n)^2}{2d_{02}r}.$$

Since $\a = \angle x_2x_0Q - \angle x_2x_0P$, the mean-value theorem yields some $\a' \in [\angle x_2x_0P,\angle x_2x_0Q ]$ with
$$\a = \frac{\cos(\angle x_2x_0Q) - \cos(\angle x_2x_0P)}{\sin(\a')} = \ell_n \frac{2d_{01} + \ell_n}{2d_{02}r\sin(\a')}.$$

To finish the proof, note that denominator is of order at least $n^{-1/4}$ by the assumptions at the beginning of proof.
\enp

We conclude the appendix with a final elementary geometric volume bound.
\bel
\label{int_min2_lem}
If $\e > 0$ is sufficiently small, then
\beals
\big|\{(x_0, x_1, x_2) \in \wl^3:\, d_{01}< d_{012} \le d_{01} + n^{-2 + \e}\}\big| \in O(n^{-1/4}).
\end{align*}
\enl

\bep
To ease notation, we put $\ell_n := n^{-2 + \e}$. Arguing similarly as in the proof of Lemma \ref{int_max_lem}, we may leverage rotational symmetry around the axis formed by $x_1x_2$ to reduce the proof to the setting $d = 2$. Similarly, we may assume that $d_{12} \ge 2\ell_n$.

%
%
 Then, for fixed $a > d_{12}$, as observed in the proof of Lemma \ref{cresc_lem}, the location of all $x_0$ with $d_{012} = a$ is the union of the two circles with radius $a$ passing through $x_1$, $x_2$. In particular, the location of $x_0$ with $d_{01} \ge a - \ell_n$ is given by an arc in each of these circles, see Figure \ref{paral_fig} (right).
Hence, applying the co-area formula \cite[Chapter 3]{coarea} to the level sets of the function $u(x_0) := d_{012}$, we obtain that 
\beals
\big|\{(x_0, x_1, x_2) \in \wl^3:\, d_{01} \ge d_{012} - \ell_n\}\big| = \int_{\ell_n}^\infty \big|\{(x_0, x_1, x_2) \in \wl^3:\, d_{01} \ge a - \ell_n\} \cap \{d_{012} = a\}\big||\nabla u(x_0)|^{-1} \d a.
\end{align*}
Since the gradient of $u$ is bounded away from 0, it suffices to show that the length of these arcs is of order $O(\sqrt{a \ell_n})$. By construction, the angle $\a$ associated with one of these arcs satisfies
$\cos(\a /2) =  (2a - 2\ell_n)/{2a},$
so that 
$$\sin^2(\a/2) = 1 - \cos^2(\a /2) \le \frac{\ell_n}a,$$
Since $\ell_n/a \le 1/2$, we deduce the asserted
$a \a \le a \sqrt{\frac{8\ell_n}a} =\sqrt{8a\ell_n}.$
\enp

\end{document}